\documentclass[fleqn, twocolumn]{article}

\usepackage[utf8]{inputenc} 
\usepackage[T1]{fontenc} 
\usepackage{textcomp}
\usepackage{lmodern}
\usepackage{lastpage} 
\usepackage{siunitx} 
\usepackage{multirow}
\usepackage{textcomp}
\usepackage{empheq}
\usepackage{pdflscape}
\usepackage{pgfplots}
\pgfplotsset{compat=1.12}
\usepackage[framed,numbered,autolinebreaks,useliterate]{mcode}
\usepackage{media9}
\usepackage{blindtext}
\usepackage[utf8]{inputenc}
\usepackage{amsmath}
\usepackage{amssymb}
\usepackage{float}
\usepackage{graphicx}
\usepackage{setspace}
\usepackage{verbatim}
\usepackage{parskip} 
\usepackage{acronym}
\usepackage{etoolbox}
\usepackage{pgfplotstable} 
\usepackage{float} 
\usepackage{xcolor}
\usepackage{multicol}
\usepackage{enumitem}
\usepackage{etoolbox}
\usepackage{graphicx}
\usepackage{comment}
\usepackage{lipsum}
\usepackage{lineno} 
\usepackage{algorithm}
\usepackage[noend]{algpseudocode}
\usepackage[margin=25mm]{geometry}
\usepackage{amsfonts}
\usepackage{verbatim}
\usepackage{subfigure}
\usepackage{multirow}
\usepackage[bookmarksnumbered=true]{hyperref} 
\hypersetup{
     colorlinks = true,
     linkcolor = blue,
     anchorcolor = blue,
     citecolor = blue,
     filecolor = blue,
     urlcolor = blue
     }
\usepackage{geometry}
\newcommand{\ep}{\varepsilon}

\newcommand{\ord}{\mathcal{O}}
\newcommand{\nw}{\noindent}
\newcommand{\der}[2] {\frac{\partial {#1} }{\partial {#2} } }

\makeatletter
\newcommand*{\rom}[1]{\expandafter\@slowromancap\romannumeral #1@}
\makeatother
\numberwithin{equation}{section}

\providecommand{\keywords}[1]
{
  \small	
  \textbf{\textit{Keywords---}} #1
}

\title{Second order adjoint sensitivity analysis in variational data assimilation for tsunami models}

\author{N. K.-R. Kevlahan$^{1}$ and R. A. Khan$^{1*}$\\
        \small $^{1}$ Department of Mathematics and Statistics, McMaster University, Hamilton, ON, Canada\\
        ${}^*$ \small{corresponding author email: ramsha.khan@math.mcmaster.ca}
}

\date{}

\begin{document}
\twocolumn[
  \begin{@twocolumnfalse}
    \maketitle
 \begin{abstract}
We mathematically derive the sensitivity of data assimilation results for tsunami  modelling, to perturbations in the observation operator. We consider results of variational  data assimilation schemes on the one dimensional shallow water equations for (i) initial condition reconstruction, and (ii) bathymetry detection as presented in Kevlahan et al. (2019, 2020) \cite{khan_2019}\cite{kevlahan2020convergence}. We use variational methods to derive the Hessian of a cost function $\mathcal{J}$ representing the error between forecast solutions and observations. Using this Hessian representation and methods outlined by Shutyaev et al. (2017, 2018) \cite{Shutyaev_ic_17}\cite{shutyaev_bath_18}, we mathematically derive the sensitivity of arbitrary response functions to perturbations in observations for case (i) and (ii) respectively. Such analyses potentially substantiate results from earlier work, on sufficient conditions for convergence \cite{khan_2019}, and sensitivity of the propagating surface wave to errors in bathymetry reconstruction \cite{kevlahan2020convergence}. Such sensitivity analyses would illustrate whether particular elements of the observation network are more critical than others, and help minimise extraneous costs for observation collection, and efficiency of predictive models. 
\end{abstract}
\vspace{15pt}
\end{@twocolumnfalse}
]

\section{\large{Introduction}}
Data assimilation methodologies are integral to accurate climate, atmospheric and oceanic modelling. Defined as optimal integration of observed data into a mathematical forecast model to refine predictions, variational data assimilation algorithms such as 3D-VAR, 4D-VAR, and Kalman filtering techniques like the EnKF (Ensemble Kalman filter) are regularly used for numerical weather prediction and longer range forecasts of climate trends. For example, data assimilation systems are used by the ECMWF for climate reanalysis, where archived observations are ``reanalysed'', in order to create a comprehensive global dataset describing the recent history of the earth's climate, atmosphere and oceans \cite{ERA5}. Additionally, data assimilation techniques are used in tsunami forecast models where observations of surface waves are used to reconstruct missing information such as initial conditions and parameters in the system, and subsequently predict impact at coastlines \cite{nakamura06}. The configuration of the observation operator has a significant impact on the results of the data assimilation scheme, as demonstrated in Kevlahan et al. (2019) \cite{khan_2019}, (2020) \cite{kevlahan2020convergence}. Therefore we have an interest in quantifying the sensitivity of the assimilation algorithm results to perturbations in the observation operator.  In this study, we utilise methods outlined in Shutyaev et al. (2017) \cite{Shutyaev_ic_17} and (2018) \cite{shutyaev_bath_18}, to mathematically derive the sensitivity of arbitrary response functions to perturbations in observations for (i) the data assimilation scheme for initial condition reconstructed outlined in \cite{khan_2019}, and (ii) the data assimilation scheme for bathymetry detection given in \cite{kevlahan2020convergence}. We do this by analytically deriving the Hessian of some cost function minimised in the data assimilation scheme, and using its uniqueness properties to derive the expressions for the sensitivity. 

The structure of this study is as follows. In section \ref{sec1dvar} we summarise the aforementioned data assimilation schemes. Section \ref{sens_IC} gives the analytical derivation of the Hessian and subsequent sensitivity analysis for initial condition assimilation results, and section \ref{sens_bath} gives the derivation for bathymetry assimilation results. Section \ref{sec_app} gives an overview of applications of these derivations to tsunami models, and concludes with consideration of current and future work. 

\section{\large{Review of first order variational data assimilation scheme}}
\label{sec1dvar}
We provide a summary of the two data assimilation schemes on the one dimensional shallow water equations,
\begin{subequations}
\begin{align}
&\frac{\partial \eta}{\partial t} + \frac{\partial }{\partial x} \Big((1 + \eta - \beta) u \Big) = 0, \\
&\frac{\partial u}{\partial t} + \frac{\partial }{\partial x} \Big( \frac{1}{2} u^2 +  \eta \Big) = 0  ,\\
&\eta(x,0) = \ \phi (x) ,\\
& u(x, 0) = \ 0,
\end{align}
\label{swe_2d_eq}
\end{subequations}
utilising a set of measurements $y^{(o)}(t)$, representing observations of the true height perturbation $\eta(x,t)$ at positions $\{ x_j\}, j=1,...,N_{obs}$.
\subsection{\small{Initial Condition Assimilation}} \label{IC_da}
For the initial condition assimilation, we assume we do not have complete information about $\phi(x)$, and our objective is to minimise the error between the forecast solution $\eta^{(f)}(x,t)$ given some guess for $\phi$,  and the observations $y^{(o)}(t)$. We define this error in terms of a cost function, 
\begin{equation}
\mathcal{J}(\phi) = \frac {1}{2} \int_0^T  \sum_{i=1}^{M} \Big[\eta^{(f)}(x_j,t;\phi) - y_j^{(o)}(t) \Big]^2 \  dt.
\label{cost_func_ic}
\end{equation}
This is equivalent to solving
\begin{equation}
    \nabla \mathcal{J}^{L^2}(\phi^{(b)}) = 0.
\end{equation}
Where $\phi^{(b)}$ is the ``best'' guess for $\phi(x)$, and the local minimiser of \eqref{cost_func_ic}. We formulate a Lagrangian constrained by \eqref{swe_2d_eq} and some arbitrarily chosen adjoint variables (Lagrange multipliers) $(\eta^*, u^*)$ that are solutions of
\begin{subequations}
\begin{empheq}[left=\empheqlbrace]{align}
& \der {\eta^*} t + {u \der {\eta^*} x}  + \der {u^*} x  \nonumber\\
&= \ H \big( \eta^{(f)}(x,t;\phi) - y^{(o)}(t) \big)  ,  \\
&\der {u^*} t + (1+{ \eta}) \der {\eta^*} x + u{ \der {u^*} x} = \ 0 , \\
& \eta^*(x,T) = \  0, \\
&\ {u^*}(x,T) = \ 0 ,
\end{empheq}\label{adj_eq_ic}
\end{subequations}
where H is the operator taking the state variable $\eta$ to the observation space. We use the Riesz representation theorem and the Gateaûx derivative representation of $\mathcal{J}'(\phi;\eta')$ given some arbitrary perturbation $\eta'$ to derive
 \begin{equation}
\nabla^{L^2} \mathcal J(\phi) = -{\eta^*}( x,0) , 
\label{grad_J_ic}
\end{equation}
where 
\begin{align}
\mathcal{J}'(\phi;\eta') \ =& \ \langle \nabla \mathcal{J}(\phi), \eta' \rangle_{{L^2}(\Omega)} \nonumber \\
=& \ \int_{-L}^{L} \nabla^{L_2} \mathcal{J}(\phi) \ \eta' \ dx.
\label{riesz_grad_IC}
\end{align}
For a more detailed analysis of this derivation we refer the reader to \cite{khan_2019}. Hence the optimal reconstruction $\phi^{(o)}$ is where 
 \begin{equation}
\nabla^{L^2} \mathcal J(\phi^{(o)}) = -{\eta^*}( x,0) = 0.
\end{equation}
In the assimilation scheme, we use a gradient descent algorithm to iteratively find the optimal reconstruction of the initial condition $\phi^{(o)}$ given some initial guess, such that \eqref{cost_func_ic} is minimised. 

\subsection{\small{Bathymetry Assimilation}} \label{bath_da}
The data assimilation scheme to recover missing bathymetry information $\beta(x)$ in \eqref{swe_2d_eq} is similar to the IC case, except now $\phi(x)$ is known, and the cost function we aim to minimise is 
\begin{equation}
\mathcal{J}(\beta) = \frac {1}{2} \int_0^T  \sum_{i=1}^{M} \Big[\eta^{(f)}(x_j,t;\beta) - y_j^{(o)}(t) \Big]^2 \  dt.
\label{cost_func_bath}
\end{equation}
Then using a similar analysis as in section \ref{IC_da}, we find that by solving the adjoint system 
\begin{subequations}
\begin{empheq}[left=\empheqlbrace]{align}
 &\der {\eta^*} t + {u \der {\eta^*} x}  + \der {u^*} x \nonumber \\
 &= \ H \big( \eta^{(f)}(x,t;\beta) - y^{(o)}(t) \big)    &,  \\
&\der {u^*} t + (1+{ \eta} - \beta) \der {\eta^*} x + u{ \der {u^*} x} = \ 0, \\
& \eta^*(x,T) = \  0, \\
& {u^*}(x,T) = \ 0 ,
\end{empheq}
\label{adj_eq_bath}
\end{subequations}
We derive 
\begin{align}
\mathcal{J}'(\beta;\beta') &= \ \langle \nabla \mathcal{J}(\beta), \beta' \rangle_{{L^2}(\Omega)} \nonumber \\
&= \ \int_{-L}^{L} \nabla^{L_2} \mathcal{J}(\beta) \ \beta' \ dx \nonumber \\
& = \int_{-L}^{L}  \int_0^{T}   u \der{\eta^*} x \beta' \ dt \ dx, 
\end{align}
hence
 \begin{equation}
\nabla^{L^2} \mathcal J = \int_0^{T}   u \der{\eta^*} x \ dt.
\label{grad_J_bath}
\end{equation}
Then the ``best'' guess for $\beta(x)$, defined as $\beta^{(o)}$ is where 
 \begin{equation}
\nabla^{L^2} \mathcal J(\beta^{(o)}) = \int_0^{T}   u \der{\eta^*} x \ dt = 0.
\end{equation}
A detailed derivation of this result can be found in \cite{kevlahan2020convergence}.

\section{\large{Initial Condition Assimilation \\ Sensitivity Analysis}} \label{sens_IC}

In order to derive the sensitivity of some arbitrary response function using the methods outlined in \cite{Shutyaev_ic_17} and \cite{shutyaev_bath_18}, we need to formulate expressions for the Hessian of the cost functions \eqref{cost_func_ic} and \eqref{cost_func_bath} formulated in section \ref{sec1dvar}. Shutyaev et al. (2017), (2018) give a general method that utilises properties of the Hessian , however does not provide a derivation, leaving it up to us to extend the adjoint analysis used to find the gradient of $\mathcal{J}$, to find the Hessian of $\mathcal{J}$. 

While works such as Le Dimet et al. (1992) \cite{ld_92} provide a derivation of the Hessian vector product for initial condition assimilation,  their  derivation is for the finite dimensional case, and assumes a vector form for both the state variables and the control variable. In our case the derivation of our first order adjoint is for the infinite dimensional case in the space $L^2(\Omega)$ over some domain $\Omega$, where we used the $L^2$ inner product and the Riesz respresentation theorem to extract our gradient $\nabla^{L^2}\mathcal{J}(\phi)$. For that reason, a derivation of the Hessian in the same functional space as for our first order adjoint is appropriate, and in doing so we aim to extract the ``Gâteaux Hessian'' for Hilbert spaces, defined thus: 

\textit{If $f$ is twice Gâteaux differentiable at $x$, we can identify $D^2f(x)$  with the operator $\nabla^2f(x) \in \mathcal{B}(\mathcal{H})$ in the sense that 
\begin{align}
    &(\forall y \in \mathcal{H})(\forall z \in \mathcal{H}), \nonumber \\  
    &\big(D^2f(x)y\big)z = \left \langle z, \nabla^2 f(x)y \right \rangle_{\mathcal{B}(\mathcal{H})},
\end{align}
where we call $\nabla^2f(x)$ the (Gâteaux) Hessian of $f$ at $x$, $\mathcal{B}(\mathcal{H})$ is the space of continuous linear functionals in $\mathcal{H}$, and $Df(x)y$ is the Gâteaux derivative of $f$ in the direction $y$.}

In the remainder of this section we present the derivation of the Hessian of $\mathcal{J}(\phi)$ and subsequently the sensitivity analysis for initial condition assimilation, using methods outlined in \cite{Shutyaev_ic_17}. Parallel results for bathymetry assimilation are presented in section \ref{sens_bath} following the methods given in \cite{shutyaev_bath_18}.

\subsection{\small{Hessian of $\mathcal{J}(\phi)$ for Initial condition \\ Assimilation}}
\nw To derive a form for the Hessian in a Hilbert space, we use the fact that we derived the following form for the Gâteaux derivative of $\mathcal{J}$ with respect to the initial condition and some perturbation direction $\eta'$ in \cite{khan_2019}:
\begin{equation}
    \mathcal{J}'(\phi;\eta') = - \int_{-L}^{L} \eta^*(x,0; \phi)\eta'\  dx.
\end{equation}
\nw Then if we consider a second perturbation of $\mathcal{J}'(\phi;\eta')$, $\eta''$ where we have $\phi \rightarrow \phi+ \ep \eta''$, the second order Gâteaux derivative of $\mathcal{J}$ can be expressed as 
\begin{align*}
    \mathcal{J}''(\phi; \eta';\eta'') =& \lim_{ \ep \to 0} \frac{\mathcal{J}'(\phi + \ep \eta'';\eta')- \mathcal{J}'(\phi;\eta')}{\ep} \\
    =&  \frac{d}{d \ep}
    \mathcal{J}'(\phi + \ep \eta'';\eta')\Bigr|_{\substack{\ep=0}}\\
    =& \frac{d}{d \ep} \Big\{ - \int_{-L}^{L} \eta^*(x,0;\phi + \ep \eta'')\eta'\  dx \Big\}.
\end{align*}
\nw We consider a regular perturbation expansion of $\eta^*(x,0;\phi + \ep \eta'')$, approximating it by the series $f_0 + f_1 \ep + \ord(\ep^2)$. We can see this is equivalent to a taylor expansion about $\ep = 0$. Then
\begin{align*}
  \frac{d}{d \ep} \mathcal{J}'(\phi+ \ep \eta'';\eta')\Bigr|_{\substack{\ep=0}} \hspace{-10pt} =& - \frac{d}{d \ep}  \int_{-L}^{L} (f_0 + f_1 \ep) \eta' \ dx \Bigr|_{\substack{\ep=0}} \\
  =& - \int_{-L}^{L} f_1 \eta' \ dx.
\end{align*}
\nw In order to find the term $f_1$, we recognise that this is equivalent to the coefficient of the linear term in the Taylor approximation, and  $f_1 = \frac{d}{d \ep} \eta^*(x,0;\phi + \ep \eta'')\Bigr|_{\substack{\ep=0}} $.

To find $f_1$, we need the adjoint variable $\eta^*$ at $t=0$ given the perturbed initial condition $\phi + \ep \eta''$. Substituting this back into our forward system \eqref{swe_2d_eq}and adjoint system \eqref{adj_eq_ic}, and gathering terms of $\ord(\ep)$ will give us our perturbed SWE and second order adjoint (SOA) system. We assume this perturbation in the initial condition brings about the following perturbations to our state and adjoint variables;
\begin{itemize}
    \item $u \rightarrow  u + \ep \hat{u} $ and $\eta \rightarrow \eta + \ep\hat{\eta}$ for the shallow water system \eqref{swe_2d_eq}.
    \item $u^* \rightarrow u^* + \ep \bar{u}$  and  $\eta^* \rightarrow \eta^* + \ep \bar{\eta}$ for the adjoint system \eqref{adj_eq_ic}.
    \end{itemize}
\nw The resulting perturbed model for the SWE is
\begin{subequations}
\begin{empheq}[left=\empheqlbrace]{align}
&\der {\hat{\eta}} t + \der {} x {\big(u\hat{\eta}\big)} + \der {} x {\big((\eta+1)\hat{u}\big)} = 0,  \\
&  \der {\hat{u}} t  + \der {\big(\hat{u} u\big)} x + \der {\hat{\eta}} x = 0,  \\
&  \hat{\eta}(x,0) = {\eta}'', \\
&  \hat{u}(x,0) = 0,
\end{empheq} \label{pert_model_H}
\end{subequations}
\nw and the second order adjoint (SOA) model is 
\begin{subequations}
\begin{empheq}[left=\empheqlbrace]{align}
 &  \der {\bar{\eta}} t + \hat{u} \der {\eta^*} x + u \der {\bar{\eta}} x + \der {\bar{u}} x = -H \hat{\eta}, & \\
 & \der {\bar{u}} t + ( \eta +1) \der {\bar{\eta}} x  + \hat{\eta}  \der {\eta^*} x \\ \nonumber
 & \hspace{14pt}+  u \der {\bar{u}} x + \hat{u} \der {u^*} x = 0, & \\
&  \bar{\eta}(x,T) = 0, \\
& \bar{u}(x,T) = 0.  & 
\end{empheq} 
\end{subequations} \label{soa_model_H}
 Going back to the definition of $f_1$, we can rewrite this as
\begin{align*}
f_1 =&  \frac{d}{d \ep} \Big( \eta^*(x,0; \phi) + \ep \bar{\eta}(x,0; \eta'') \Big)\Bigr|_{\substack{\ep=0}} \\
=& \bar{\eta}(x,0; \eta''),
\end{align*}
And thus we have
\begin{align}
 \mathcal{J}''(\eta_0; \eta';\eta'') =&   - \int_{-L}^{L}\bar{\eta}(x,0;\eta'')\eta'\  dx  = \\
 &\left \langle - \bar{\eta}(x,0;\eta''), \ \eta'  \right \rangle_{L^2(\Omega)}
\end{align}
and by the definition of the Gâteaux Hessian, we get $\nabla^2 \mathcal{J}(\eta_0;\eta'') = - \bar{\eta}(x,0;\eta'') $.

Extending the kappa test analysis outlined in \cite{khan_2019} to ensure our derivation is correct, we define
\begin{equation}
    \kappa(\ep) = \lim_{ \ep \to 0} \frac{1}{\ep} \frac{\mathcal{J}'(\phi + \ep \eta'';\eta')- \mathcal{J}'(\phi;\eta')}{- \int_{-L}^{L}\bar{\eta}(x,0;\eta'')\eta'\  dx },
\end{equation}
\nw where 
\begin{align}
    &\ \mathcal{J}'(\eta_0 + \ep \eta'';\eta') \\
 = & \lim_{ \tau \to 0}  \frac{\mathcal{J}'(\phi + \ep \eta'' + \tau \eta')- \mathcal{J}'(\phi + \ep \eta'')}{\tau},
\end{align}
\nw and
\begin{equation}
 \mathcal{J}'(\phi;\eta') =    \lim_{ \tau \to 0}  \frac{\mathcal{J}'(\phi + \tau \eta')- \mathcal{J}'(\phi)}{\tau}.
\end{equation}
\nw If the derivations are correct, $\kappa(\ep)$ should converge to $1$ as $\ep \rightarrow 0$.

A key observation here is that we have only derived the action of the Hessian on some perturbation $\eta''$, whereas in section \ref{IC_da} we were able to derive the gradient of $\mathcal{J}(\phi)$ for any arbitrary perturbation $\eta'$. It does not seem possible to proceed as before and find the Hessian for any arbitrary $\eta''$ using variational methods. Indeed, \cite{ld_92} were also limited to derivation of a ``Hessian vector product'' in the finite dimensional case. Despite this, the current derivation is satisfactory for the following sensitivity analysis. 

\subsection{\small{Sensitivity Analysis for Initial Condition Assimilation}}

Define the optimality system as the successive solution of 
\begin{subequations}
\begin{empheq}[left=\empheqlbrace]{align}
&\frac{\partial \eta}{\partial t} + \frac{\partial }{\partial x} \Big((1 + \eta ) u \Big) = 0, \\
&\frac{\partial u}{\partial t} + \frac{\partial }{\partial x} \Big( \frac{1}{2} u^2 +  \eta \Big) = 0  ,\\
&\eta(x,0) = \ \phi (x) ,\\ 
& u(x, 0) = \ 0. 
\end{empheq}
\label{swe_2d_ic_opt}
\end{subequations}
\begin{subequations}
\hspace{-20pt}
\begin{empheq}[left=\empheqlbrace]{align}
& \der {\eta^*} t + {u \der {\eta^*} x}  + \der {u^*} x = H\big( \eta - m\big)   ,  \\
&\der {u^*} t + (1+{ \eta}) \der {\eta^*} x + u{ \der {u^*} x} = \ 0 , \\
& \eta^*(x,T) = \  0, \\
\label{adj_ic}
& {u^*}(x,T) = \ 0 ,
\end{empheq}
\label{adj_eq_ic_opt}
\end{subequations}
Where $\phi(x)$ is the optimal reconstruction of the initial condition $\eta(x,0)$ giving 
 \begin{equation}
- \eta^*(x,0) = 0,
\label{grad_J}
\end{equation}
and $m(t)$ are the observations taken at $\{x_j\}$ at for $j = 1,...,N_{obs}$.
Let us consider some arbitrary response function
\begin{equation}
    \mathcal{G}(\eta, u, \phi).
\end{equation}
Then by the chain rule, the sensitivity of $\mathcal{G}$ to perturbations in the observations $m$ can be defined as 
\begin{equation}
  \frac{d\mathcal{G}}{d m} = \frac{d\mathcal{G}}{d \eta}\frac{d\eta}{d m} + \frac{d\mathcal{G}}{d u}\frac{d u}{d m} + \frac{d\mathcal{G}}{d \phi}\frac{d\phi}{d m}  .
\end{equation}
Let us now consider a perturbation in the observations, such that given $m \rightarrow m + \hat{m}$ gives
\begin{itemize}
    \item $u \rightarrow  u + \hat{u} $ and $\eta \rightarrow \eta + \hat{\eta}$ for the shallow water system \eqref{swe_2d_ic_opt}.
    \item $u^* \rightarrow u^* + {\tilde{u}}^*$  and  $\eta^* \rightarrow \eta^* + {\tilde{\eta}}^*$ for the adjoint system \eqref{adj_eq_ic_opt}.
     \item $\phi \rightarrow \phi + \hat{\phi}$ for the optimal initial condition.
    \end{itemize}
Our perturbed system becomes
\begin{subequations}
\begin{empheq}[left=\empheqlbrace]{align}
&\der {\hat{\eta}} t + \der {} x {\big(u\hat{\eta}\big)} + \der {} x {\big((\eta+1)\hat{u}\big)} = 0,  \\
&  \der {\hat{u}} t  + \der {\big(\hat{u} u\big)} x + \der {\hat{\eta}} x = 0,  \\
&  \hat{\eta}(x,0) = \hat{\phi}(x), \\
&  \hat{u}(x,0) = 0,
\end{empheq}\label{[pert_m_swe}
\end{subequations}
\begin{subequations}
\begin{empheq}[left=\empheqlbrace]{align}
 &  \der {{\tilde{\eta}}^*} t + \hat{u} \der {\eta^*} x + u \der {{\tilde{\eta}}^*} x + \der {{\tilde{u}}^*} x \\ \nonumber
 & \hspace{14pt} = H (\hat{m} - \hat{\eta}(x_i,t ;\hat{\phi}), & \\
 & \der {{\tilde{u}}^*} t + ( \eta +1) \der {{\tilde{\eta}}^*} x  + \hat{\eta}  \der {\eta^*} x \\ \nonumber
 & \hspace{14pt}+  u \der {{\tilde{u}}^*} x + \hat{u} \der {u^*} x = 0, & \\
&  {\tilde{\eta}}^*(x,T) = 0, \\
& {\tilde{u}}^*(x,T) = 0.  &
\end{empheq}
\label{pert_m_adj}
\end{subequations}
\begin{equation}
 -{ \hat{\eta}}^*(x,0) = 0.
 \label{pert_grad}
\end{equation}
Then we can say 
\begin{equation}
\Bigg\langle \frac{d\mathcal{G}}{d m}, \hat{m}\Bigg\rangle_{Y_{obs}}\hspace{-15pt} = \Bigg\langle \der {\mathcal{G}}{\eta}, \hat{\eta}\Bigg\rangle_{Y} + \Bigg\langle \der{\mathcal{G}}{u}, \hat{u} \Bigg\rangle_{Y} 
+\Bigg\langle \der{\mathcal{G}}{\phi}, \hat{\phi} \Bigg\rangle_{Y_p}
 \vspace{-1pt}
\label{chain_rule}
\end{equation}
where 
\begin{align}
   Y_{obs} = L^2\Big([-L,L] \times[0,T] \Big), \\
   Y = L^2\Big([-L,L] \times[0,T] \Big), \\
   Y_{p} \in L^2\Big([-L,L] \Big),
\end{align}
are the observation space, state space, and initial condition space respectively. 
Let us introduce some adjoint variables $P_i$, $i = 1,...,5$, where $P_i \in Y$, $i = 1,...,4$, and$P_5 \in Y_{obs}$. Then if we take the inner product of $P_1$ and $P_2$ with the systems \eqref{[pert_m_swe}, inner product of $P_3$ and $P_4$ with \eqref{pert_m_adj}, and $P_5$ with \eqref{pert_grad}, we get the following duality relation
\begin{align}
    &0 =  \nonumber\\ 
    &\ \int_{0} ^{T}\int_{-L}^L \Bigg\{ P_1 \Bigg[  \der {\hat{\eta}} t + \der {} x {\big(u\hat{\eta}\big)} + \der {} x {\big((\eta+1)\hat{u}\big)}  \Big ] \nonumber\\ 
         &+  P_2 \Big [ \der {\hat{u}} t  + \der {\big(\hat{u} u\big)} x + \der {\hat{\eta}} x    \Big ]  \nonumber\\ 
         &+  P_3 \Big [   \der {{\tilde{\eta}}^*} t + \hat{u} \der {\eta^*} x + u \der {{\tilde{\eta}}^*} x + \der {{\tilde{u}}^*} x  - H (\hat{m} - \hat{\eta}) \Big ] \nonumber  \\ 
         &+  P_4 \Big [  \der {{\tilde{u}}^*} t + ( \eta +1) \der {{\tilde{\eta}}^*} x  + \hat{\eta}  \der {\eta^*} x \nonumber \\ 
         &+  u \der {{\tilde{u}}^*} x + \hat{u} \der {u^*} x\Big ]   \Bigg\} dt dx \nonumber\\ 
         &+ \int_{-L}^L  P_5 \Big [  -{ \hat{\eta}}^*(x,0) \Big ] dx. \label{adj_ibp}
\end{align}

integrating the double integral in \eqref{adj_ibp} by parts in space and time, we are able to transfer the derivatives onto the adjoint variables $P_i$, $i = 1,...,5$ instead of on $\hat{u}, \hat{\eta}, {\tilde{u}^*}, {\tilde{u}^*}$. Since the choice of adjoint variables is arbitrary, we pick the following systems for $P_i$,
\begin{subequations}
\begin{empheq}[left=\empheqlbrace]{align}
&\der {P_3} t + \der {} x {\big(uP_3\big)} + \der {} x {\big((\eta+1)P_4\big)} = 0,  \\
&  \der {P_4} t  + \der {\big(u P_4 \big)} x + \der {P_3} x = 0,  \\
&  P_3(x,0) = -P_5 \\
&  P_4(x,0) = 0
\end{empheq} \label{p3_4}
\end{subequations}
\begin{subequations}
\begin{empheq}[left=\empheqlbrace]{align}
 &  \der {P_1} t + P_4 \der {\eta^*} x + u \der {P_1} x \nonumber \\
 &\hspace{20pt}+ \der {P_2} x -H P_3 = - \der{\mathcal{G}}{\eta}, & \\
 & \der {P_2} t + ( \eta +1) \der {P_1} x  + P_3  \der {\eta^*} x \\ \nonumber
 & \hspace{14pt}+  u \der {P_2} x + P_4 \der {u^*} x = - \der{\mathcal{G}}{u}, & \\
&  P_1(x,T) =  P_5  - \der{\mathcal{G}}{\phi}, \\
& P_2(x,T) = 0.  & 
\end{empheq} \label{p1_2}
\end{subequations}
We subsequently get rid of $P_5$ by utilising the fact that $P_1(x,T) =  P_5  - \der{\mathcal{G}}{\phi}$, and $P_3(x,0) = -P_5$, and define $P_3(x,0) = \nu$, where the auxiliary variable $\nu$ is defined as 
\begin{equation}
    \nu = \der{\mathcal{G}}{\phi} - P_1(x,0).
\end{equation}
 Then subsequently, as a result of integration by parts and the choice of systems for $P_i$, \eqref{adj_ibp} reduces to
 \begin{equation}
 \Bigg\langle HP_3, \hat{m}\Bigg\rangle_{Y_{obs}}\hspace{-15pt} = \Bigg\langle \der {\mathcal{G}}{\eta}, \hat{\eta}\Bigg\rangle_{Y} \hspace{-8pt} + \Bigg\langle \der{\mathcal{G}}{u}, \hat{u} \Bigg\rangle_{Y} \hspace{-8pt} + \Bigg\langle \der{\mathcal{G}}{\phi}, \hat{\phi} \Bigg\rangle_{Y_p}\hspace{-8pt} .
 \label{riesz_2}
 \end{equation}
By the Riesz representation theorem and equivalence of inner products in \eqref{chain_rule} and \eqref{riesz_2}, we can define the sensitivity of the response function $\mathcal{G}(\eta, u, \phi)$ to perturbations in the observations $m$ as 
 \begin{equation}
     \der {\mathcal{G}}m = HP_3(x,t). \label{sensitivity}
 \end{equation}
Finding $P_3$ requires solving the systems \eqref{p3_4} (with $P_3(x,0) = \nu$), and \eqref{p1_2}. This is a coupled system of four variables, with two initial time conditions and two final time conditions, making it challenging to solve. However, we see that \eqref{p3_4} is equivalent to the perturbed system for the Hessian $\nabla^2 \mathcal{J}(\phi)$ \eqref{pert_model_H}, and \eqref{p1_2} is equivalent to the second order adjoint (SOA) system \eqref{soa_model_H} with forcing term $\Big(- \der{\mathcal{G}} {\eta}, - \der{\mathcal{G}} {u}  \Big)^T$.
Shutyaev et al. (2017) illustrate that the solutions to the adjoint systems \eqref{p3_4} and \eqref{p1_2} are then equivalent to solving 
\begin{equation}
    \mathcal{H}\nu = \mathcal{F},\label{op_eq_ic}
\end{equation} 
where $\mathcal{F}$ is defined as 
\begin{equation}
    \mathcal{F} = \der{\mathcal{G}}{\phi} + \psi(x,0),
\end{equation}
and $\psi$ is the solution of the forced first order adjoint system 
\begin{subequations}
\hspace{-20pt}
\begin{empheq}[left=\empheqlbrace]{align}
& \der {\psi} t + {u \der {\psi} x}  + \der {\varphi} x = \der{\mathcal{G}}{\eta}  ,  \\
&\der {\varphi} t + (1+{ \eta}) \der {\psi} x + u{ \der {\varphi} x} = \der{\mathcal{G}}{u}  , \\
& \psi(x,T) = \  0, \\
& {\varphi}(x,T) = \ 0 .
\end{empheq}\label{foa_forced}
\end{subequations}

The significance of recognising the systems \eqref{p1_2} and \eqref{p3_4} as the Hessian with external forcing, is that if we assume $\mathcal{H}$ is positive definite, the operator equation \eqref{op_eq_ic} is correctly and everywhere solvable in $Y$ \cite{Shutyaev_ic_17}. Hence for every $\mathcal{F}$, we can find a unique $\nu$ such that \eqref{op_eq_ic} holds. Then the sensitivity $\der {\mathcal{G}}m$ can be found by the following steps:
\begin{enumerate}
    \item Define $\mathcal{F} = \der{\mathcal{G}}{\phi} + \psi(x,0)$, where $\psi$ is the solution of \eqref{foa_forced}.
    \item Solve $\mathcal{H}\nu = \mathcal{F}$ for $\nu$.
    \item Solve the system \eqref{p3_4} using $P_3(x,0) = \nu$ to find $P_3(x,t)$.
    \item Define $\der {\mathcal{G}}m = HP_3(x,t)$, where $H$ is the operator mapping the $\eta$ from state space $Y$ to the observation space $Y_{obs}$ .
\end{enumerate}

\section{\large{Bathymetry assimilation sensitivity analysis}}\label{sens_bath}

\subsection{\small{Hessian of $\mathcal{J} (\beta)$}}

The derivation of the Hessian for $\mathcal{J}(\beta)$ and subsequent sensitivity analysis for the bathymetry assimilation scheme described in  section \ref{bath_da} is parallel to the initial condition case. We verify this below. 

We know from section \ref{bath_da} the Gâteaux derivative of $\mathcal{J}$ with respect to the bathymetry and some perturbation direction $\beta'$ is
\begin{equation}
    \mathcal{J}'(\beta;\beta') = - \int_{-L}^{L} \Bigg(\int_0^T u \der{\eta^*} x  \ dt\Bigg) \beta'\  dx.
\end{equation}
Consider a second perturbation of $\mathcal{J}'(\beta;\beta')$, $\hat{\beta}$ where we have $\beta \rightarrow \beta + \ep \hat{\beta}$. Then the second order Gâteaux derivative of $\mathcal{J}$ is
\begin{align*}
    & \mathcal{J}''(\beta; \beta';\hat{\beta}) \\
    =&  \frac{d}{d \ep}
    \mathcal{J}'(\beta + \ep \hat{\beta};\beta')\Bigr|_{\substack{\ep=0}}\\
    =& \frac{d}{d \ep} \Big\{  \int_{-L}^{L} \int_0^T \hspace{-5pt} u(x,t;\beta+\ep\hat{\beta}) \der{\eta^*}x (x,t;\beta+ \ep \hat{\beta})  \ dt dx \Big\}\Bigr|_{\substack{\ep=0}}.
\end{align*}
Considering a regular perturbation expansion of the integrand as before, we approximate it by the series $f_0 + f_1 \ep + \ord(\ep^2)$. Then
\vspace{-5pt}
\begin{align*}
  & \frac{d}{d \ep} \mathcal{J}'(\beta + \ep \hat{\beta};\beta')\Bigr|_{\substack{\ep=0}} \\
  = -&  \frac{d}{d \ep}  \int_{-L}^{L} \Big(\int_0^T(f_0 + f_1 \ep) \  dt \Big) \beta' \ dx \Bigr|_{\substack{\ep=0}} \\
  = - &\int_{-L}^{L} \Big(\int_0^T f_1 \  dt \Big) \beta' \ dx.
\end{align*}

 We have $f_1 = \frac{d}{d \ep} u(x,t;\beta+\ep\hat{\beta}) \der{} x \eta^*(x,t;\beta+ \ep \hat{\beta})\Bigr|_{\substack{\ep=0}} $, and require $u$ and  adjoint variable $\eta^*$ given the perturbed bathymetry $\beta + \ep \hat{\beta}$. To find the resulting forward and adjoint system given the perturbation, we assume this perturbation in the initial condition brings about the following perturbations to our state and adjoint variables;
\begin{itemize}
    \item $u \rightarrow  u + \hat{u} $ and $\eta \rightarrow \eta + \hat{\eta}$ for the shallow water system.
    \item $u^* \rightarrow u^* + \bar{u}$  and  $\eta^* \rightarrow \eta^* + \bar{\eta}$ for the adjoint system.
    \end{itemize}
The resulting perturbed model for the state variables $\hat{u}, \hat{\eta}$ is
\begin{subequations}
\begin{empheq}[left=\empheqlbrace]{align}
 \der {\hat{\eta}} t &+ \der{}x \big( (1+ \eta - \beta)\hat{u} \big) \\ \nonumber
& + \der {\big(u\hat{\eta}\big)} x  
   - \der {\big(\hat{\beta} u\big)} x = 0,  \\
 \der {\hat{u}} t  &+  \der {\big(\hat{u} u\big)} x + \der {\hat{\eta}} x = 0, \\
  \hat{\eta}(&x,0) = 0, \\ 
   \hat{u}(&x,0) = 0.
\end{empheq}
\end{subequations} \label{pert_model_H_bath}

\nw and the second order adjoint (SOA) model is 
\begin{subequations}
\begin{empheq}[left=\empheqlbrace]{align}
   \der {\bar{\eta}} t &+ \hat{u} \der {\eta^*} x + u \der {\bar{\eta}} x + \der {\bar{u}} x = - \hat{\eta}(x_i,t ;\hat{\beta}),  \\
   \der {\bar{u}} t  &+ ( 1 + \eta - \beta) \der {\bar{\eta}} x + (\hat{\eta} - \hat{\beta} ) \der {\eta^*} x \\ \nonumber
& + u \der {\bar{u}} x +  \hat{u} \der {u^*} x  = 0, \\
&  \hspace{-10pt}\bar{\eta}(x,T) = 0, \\
& \hspace{-10pt} \bar{u}(x,T) = 0. 
\end{empheq}
\end{subequations} \label{adj_model_H_bath}
Giving us
\begin{equation}
  \nabla^2 \mathcal{J}(\beta;\hat{\beta}) = \int_0^T \big(  \hat{u} \der{{\eta^*}} x + u \der{\bar{\eta}} x \big)  dt.   
\end{equation}
We define the hessian $\mathcal{H}$  acting on the perturbation $\hat{\beta}$ as the successive solution of the perturbed and SOA models such that
\begin{equation}
  \mathcal{H}\hat{\beta} =  \int_0^T \big(  \hat{u} \der{{\eta^*}} x + u \der{\bar{\eta}} x \big)  dt . 
\end{equation}
This derivation can be verified using the two different forms of $\nabla^2 \mathcal{J}(\beta; \beta';\hat{\beta})$:

\begin{equation}
    \kappa(\ep) = \lim_{ \ep \to 0}  \frac{1}{\ep} \frac{\mathcal{J}'(\beta + \ep \hat{\beta};\beta')- \mathcal{J}'(\beta;\eta')}{- \int_{-L}^{L}  \Big( \int_0^T \big(  \hat{u} \der{{\eta^*}} x + u \der{\bar{\eta}} x \big)  dt \Big) \beta'\  dx },
\end{equation}
As before, if the derivations are correct, $\kappa(\ep)$ should converge to $1$ as $\ep \rightarrow 0$.

\subsection{\small{Sensitivity analysis for bathymetry assimilation}}
Define the optimality system for bathymetry assimilation as the successive solution of
\begin{subequations}
\begin{empheq}[left=\empheqlbrace]{align}
&\frac{\partial \eta}{\partial t} + \frac{\partial }{\partial x} \Big((1 + \eta ) u \Big) = 0, \\
&\frac{\partial u}{\partial t} + \frac{\partial }{\partial x} \Big( \frac{1}{2} u^2 + \eta \Big) = 0  ,\\
&\eta(x,0) = \ \eta_0(x) ,\\ 
& u(x, 0) = \ 0. 
\end{empheq}
\label{swe_2d_bath_opt}
\end{subequations}
\begin{subequations}
\hspace{-20pt}
\begin{empheq}[left=\empheqlbrace]{align}
 &\der {\eta^*} t + {u \der {\eta^*} x}  + \der {u^*} x = \ H\big( \eta - m \big) ,  \\
&\der {u^*} t + (1+{ \eta} - \beta) \der {\eta^*} x + u{ \der {u^*} x} = \ 0 , \\
&\eta^*(x,T) = \  0, \\
\label{adj_ic}
& \ {u^*}(x,T) = \ 0 ,
\end{empheq}
\label{adj_eq_bath_opt}
\end{subequations}
Where $\lambda(x)$ is the optimal reconstruction of the bathymetry $\beta(x)$ giving :
 \begin{equation}
 \int_0^{T}   u \der{\eta^*} x \ dt = 0.
\label{grad_J_bath}
\end{equation}

Let us consider some arbitrary response function
\begin{equation}
    \mathcal{G}(\eta, u, \lambda).
\end{equation}
Then by the chain rule, the sensitivity of $\mathcal{G}$ to perturbations in the observations $m$ can be defined as 
\begin{equation}
  \frac{d\mathcal{G}}{d m} = \frac{d\mathcal{G}}{d \eta}\frac{d\eta}{d m} + \frac{d\mathcal{G}}{d u}\frac{d u}{d m} + \frac{d\mathcal{G}}{d \lambda}\frac{d\lambda}{d m}  .
\end{equation}
Give some perturbation of the observations such that given $m \rightarrow m + \hat{m}$ we have
\begin{itemize}
    \item $u \rightarrow  u + \hat{u} $ and $\eta \rightarrow \eta + \hat{\eta}$ for the shallow water system \eqref{swe_2d_bath_opt}.
    \item $u^* \rightarrow u^* + {\tilde{u}}^*$  and  $\eta^* \rightarrow \eta^* + {\tilde{\eta}}^*$ for the adjoint system \eqref{adj_eq_bath_opt}.
     \item $\lambda \rightarrow \phi + \hat{\lambda}$ for the optimal initial condition.
    \end{itemize}
Our perturbed system becomes
\begin{subequations}
\begin{empheq}[left=\empheqlbrace]{align}
&\der {\hat{\eta}} t + \der {} x {\big(u\hat{\eta}\big)} + \der {} x {\big((\eta+1- \lambda)\hat{u}\big)} \nonumber\\ 
&- \der {} x {\big(\hat{\lambda} u\big)} - \der {} x {\big(\lambda\hat{u}\big)}= 0,  \\
&  \der {\hat{u}} t  + \der {\big(\hat{u} u\big)} x + \der {\hat{\eta}} x = 0,  \\
&  \hat{\eta}(x,0) = \hat{\phi}(x), \\
&  \hat{u}(x,0) = 0,
\end{empheq}\label{[pert_m_swe_bath}
\end{subequations}
\begin{subequations}
\begin{empheq}[left=\empheqlbrace]{align}
 &  \der {{\tilde{\eta}}^*} t + \hat{u} \der {\eta^*} x + u \der {{\tilde{\eta}}^*} x + \der {{\tilde{u}}^*} x \\ \nonumber
 & \hspace{18pt} = H (\hat{m} - \hat{\eta}(x_i,t ;\hat{\lambda}), & \\
 & \der {{\tilde{u}}^*} t + ( \eta +1 - \lambda) \der {{\tilde{\eta}}^*} x  + \hat{\eta}  \der {\eta^*} x \\ \nonumber
 & +  u \der {{\tilde{u}}^*} x + \hat{u} \der {u^*} x   - \hat{\lambda} \der {\eta^*} x  = 0, & \\
&  {\tilde{\eta}}^*(x,T) = 0, {\tilde{u}}^*(x,T) = 0.  &
\end{empheq}
\label{pert_m_adj_bath}
\end{subequations}
\begin{equation}
 \int_0^T \big(  \hat{u} \der{{\eta^*}} x + u \der{{\tilde{\eta}}^*} x \big)  dt . 
 \label{pert_grad_bath}
\end{equation}
Then we can say 
\begin{equation}
\Bigg\langle \frac{d\mathcal{G}}{d m}, \hat{m}\Bigg\rangle_{Y_{obs}}\hspace{-15pt} = \Bigg\langle \der {\mathcal{G}}{\eta}, \hat{\eta}\Bigg\rangle_{Y} + \Bigg\langle \der{\mathcal{G}}{u}, \hat{u} \Bigg\rangle_{Y} 
+\Bigg\langle \der{\mathcal{G}}{\lambda}, \hat{\lambda} \Bigg\rangle_{Y_p}
 \vspace{-1pt}
\label{chain_rule_bath}
\end{equation}
where 
\begin{align}
   Y_{obs} = L^2\Big([-L,L] \times[0,T] \Big) ,\\
   Y = L^2\Big([-L,L] \times[0,T] \Big), \\
   Y_{p} = L^2\Big([-L,L] \Big)
\end{align}
are as before in section \ref{sens_IC}. Let us introduce some adjoint variables $P_i$, $i = 1,...,5$, where $P_i \in Y$, $i = 1,...,4$, and $P_5 \in Y_{obs}$. Then if we take the inner product of $P_1$ and $P_2$ with the systems \eqref{[pert_m_swe_bath}, inner product of $P_3$ and $P_4$ with \eqref{pert_m_adj_bath}, and $P_5$ with \eqref{pert_grad_bath}, we get the following duality relation
\begin{align}
    &0 =  \nonumber\\ 
    &\ \int_{0} ^{T}\int_{-L}^L \Bigg\{ P_1 \Bigg[  \der {\hat{\eta}} t + \der {} x {\big(u\hat{\eta}\big)} + \der {} x {\big((\eta+1- \lambda)\hat{u}\big)} \nonumber\\ 
&- \der {} x {\big(\hat{\lambda} u\big)} - \der {} x {\big(\lambda\hat{u}\big)}  \Big ] \nonumber\\ 
         &+  P_2 \Big [  \der {\hat{u}} t  + \der {\big(\hat{u} u\big)} x + \der {\hat{\eta}} x    \Big ]  \nonumber\\ 
         &+  P_3 \Big [   \der {{\tilde{\eta}}^*} t + \hat{u} \der {\eta^*} x + u \der {{\tilde{\eta}}^*} x + \der {{\tilde{u}}^*} x \\ \nonumber
 & - H (\hat{m} - \hat{\eta}(x_i,t ;\hat{\lambda}) \Big ] \nonumber  \\ 
         &+  P_4 \Big [ \der {{\tilde{u}}^*} t + ( \eta +1 - \lambda) \der {{\tilde{\eta}}^*} x  + \hat{\eta}  \der {\eta^*} x \\ \nonumber
 & +  u \der {{\tilde{u}}^*} x + \hat{u} \der {u^*} x   - \hat{\lambda} \der {\eta^*} x\Big ]    \nonumber\\ 
         &+  P_5 \Big [  \hat{u} \der{{\eta^*}} x + u \der{{\tilde{\eta}}^*} x  \Big ] \Bigg\} dt dx. \label{adj_ibp_bath}
\end{align}
Integrating \eqref{adj_ibp_bath} by parts in space and time, we are able to transfer the derivativatives onto the adjoint variables $P_i$, $i = 1,...,5$ instead of on $\hat{u}, \hat{\eta}, {\tilde{u}^*}, {\tilde{u}^*}$. We pick the following systems for $P_i$,
\begin{subequations}
\begin{empheq}[left=\empheqlbrace]{align}
 \der {P_3} t &+ \der{}x \big( (1+ \eta - \lambda)P_4 \big) \\ \nonumber
& - \der {\big(uP_5\big)} x  
   - \der {\big(\lambda P_4\big)} x = 0,  \\
 \der {P_4} t  &+  \der {\big(P_4 u\big)} x + \der {P_3} x = 0, \\
  P_3(&x,0) = 0, \\
  P_4(&x,0) = 0.
\end{empheq}\label{p3_4_bath}
\end{subequations} 
\begin{subequations}
\begin{empheq}[left=\empheqlbrace]{align}
   \der {P_1} t &+ P_4 \der {\eta^*} x + u \der {P_1} x  \nonumber \\
   &+ \der {P_2} x  - HP_3 =  - \der{\mathcal{G}} {\eta},  \\
   \der {P_2} t  &+ ( 1 + \eta - \lambda) \der {P_1} x + (P_3 - P_5) \der {\eta^*} x \\ \nonumber
& + u \der {P_2} x +  P_4 \der {u^*} x  = - \der{\mathcal{G}} {u}, \\
&  P_1(x,T) = 0, \\
& P_2(x,T) = 0. 
\end{empheq}\label{p1_2_bath}
\end{subequations} 
\begin{equation}
 \int_0^T \big(  P_4 \der{{\eta^*}} x + u \der{P_1} x \big)  dt = \der{\mathcal{G}}{\lambda}. 
 \label{grad_adj_bath}
\end{equation}
Subsequently as a result of integration by parts and the choice of systems for $P_i$, \eqref{adj_ibp_bath} reduces to
 \begin{equation}
 \Bigg\langle P_3, \hat{m}\Bigg\rangle_{Y_{obs}}\hspace{-15pt} = \Bigg\langle \der {\mathcal{G}}{\eta}, \hat{\eta}\Bigg\rangle_{Y} \hspace{-8pt} + \Bigg\langle \der{\mathcal{G}}{u}, \hat{u} \Bigg\rangle_{Y} \hspace{-8pt} + \Bigg\langle \der{\mathcal{G}}{\lambda}, \hat{\lambda} \Bigg\rangle_{Y_p}\hspace{-8pt} .
 \label{riesz_2_bath}
 \end{equation}
By the Riesz representation theorem and equivalence of inner products in \eqref{chain_rule} and \eqref{riesz_2}, we  define the sensitivity of the response function $\mathcal{G}(\eta, u, \lambda)$ to perturbations in the observations $m$ as 
 \begin{equation}
     \der {\mathcal{G}}m = HP_3(x,t). \label{sensitivity_bath}
 \end{equation}
As in the initial condition case, we observe that \eqref{p3_4_bath} is equivalent to the perturbed system for the Hessian $\nabla^2 \mathcal{J}(\phi)$ \eqref{pert_model_H_bath} with $P_5 = \hat{\beta}(x)$, and \eqref{p1_2_bath} is equivalent to the second order adjoint (SOA) system \eqref{adj_model_H_bath} with forcing term $\Big(- \der{\mathcal{G}} {\eta}, - \der{\mathcal{G}} {u}  \Big)^T$. Let us replace $P_5$ with the auxiliary variable $\nu$. Then Shutyaev et al. (2018) \cite{shutyaev_bath_18} illustrate that the solutions to the adjoint systems \eqref{p3_4_bath} and \eqref{p1_2_bath} are equivalent to solving 
\begin{equation}
    \mathcal{H}\nu = \mathcal{F},
\end{equation}
where $\mathcal{F}$ is defined as 
\begin{equation}
    \mathcal{F} = \der{\mathcal{G}}{\lambda} - \int_0^T u \der{\gamma} x dt,
\end{equation}
and $\gamma$ is the solution of the forced first order adjoint system 
\begin{subequations}
\begin{empheq}[left=\empheqlbrace]{align}
&\der {\gamma} t + {u \der {\gamma} x}  + \der {\psi} x = - \der{\mathcal{G}}{\eta}, \\
&\der {\psi} t + (1+{ \eta} - \beta) \der {\gamma} x + u{ \der {\psi} x} = - \der{\mathcal{G}}{u},\\
&\gamma(x,T) = 0,\\
&{\psi}(x,T) = \ 0 .
\end{empheq}\label{forced_adj_bath}
\end{subequations}
As before, under the assumption $\mathcal{H}$ is positive definite, we can find a unique $\nu$ for every $\mathcal{F}$ such that $ \mathcal{H}\nu = \mathcal{F}$. The we find $\der {\mathcal{G}}m$ by:

\begin{enumerate}
    \item Defining  $\der{\mathcal{G}}{\lambda} - \int_0^T u \der{\gamma} x dt$, where $\gamma$ is the solution of \eqref{forced_adj_bath}.
    \item Solve $\mathcal{H}\nu = \mathcal{F}$ for $\nu$.
    \item Solve the system \eqref{p3_4} using $P_5(x) = \nu$ to find $P_3(x,t)$.
    \item Define $\der {\mathcal{G}}m = HP_3(x,t)$, where $H$ is the operator mapping the $\eta$ from state space $Y$ to the observation space $Y_{obs}$ .
\end{enumerate}

\section{\large{Applications and further considerations}}\label{sec_app}

The results of this analysis would be a significant addition to the work undertaken in Kevlahan et al. (2019) \cite{khan_2019}, and (2020) \cite{kevlahan2020convergence}. In the former, we implemented a data assimilation scheme on the shallow water equations, and derived a theorem giving sufficient conditions for convergence to the true initial condition; at least one pair of observation points must
be spaced more closely than half the effective minimum wavelength of the energy spectrum of the initial conditions. This conclusion can be further substantiated by considering the effects of a perturbation on the observations that defies this condition, where tge response function $\mathcal{G}$ represents the $L^2$ error between the exact and optimally reconstructed initial condition. 

Additionally, in the data assimilation scheme implemented for bathymetry detection in \cite{kevlahan2020convergence}, we illustrated two key results; reconstruction of the bathymetry is worse when the first surface wave measurement is taken after bathymetry has been observed, and there is low sensitivity of the surface wave to errors in the bathymetry reconstruction.  If  $\mathcal{G}$ represents the error in bathymetry reconstruction, by considering perturbations in the observations placed before and after the support of the bathymetry, we can gain further insight on the former result. The low sensitivity of the propagating surface wave  to bathymetry reconstruction error is an especially poignant result in the context of tsunami propagation. If our primary objective is effective prediction of the surface wave given some bathymetry reconstruction, then we can have relaxed criteria for convergence in the data assimilation scheme. We need further investigation to substantiate this result, and this sensitivity analysis is a very useful way to achieve this. It would illustrate whether certain elements of the observation network are more critical than others, and help minimise extraneous costs for observation collection and efficiency in our predictive models. 

Numerical implementations of this sensitivity analysis are currently being studied. Given that that the system $\mathcal{H}\nu = \mathcal{F}$ contains more unknowns than equations, we cannot use discretisation methods to solve for $\nu$, and subsequently we use root finding iterative  processes to find $\nu$ such that $\parallel \mathcal{H}\nu - \mathcal{F} \parallel_{\infty} = 0$. 

Further work would be to consider sensitivity of the data assimilation results to perturbations in parameters of the system, instead of the observations. Current methods being considered are variance based sensitivity analysis, that allows the variance of the scheme to be decomposed into fractions attributed to individual parameters \cite{sobol}. These results in combination with the results presented in this study would build on top of the work in \cite{khan_2019} and \cite{kevlahan2020convergence}, and pave the way for more effective tsunami models applicable to realistic scenarios.

\bibliographystyle{plain}
\bibliography{bibliography}
\nocite{*}

\end{document}